\documentclass{article}
\usepackage{amsthm,amsfonts,amsmath,amssymb,hyperref}
\usepackage[numbers]{natbib}

\newcommand{\Mohle}{{M\"ohle}}

\newcommand{\ed}{\stackrel{d}{=}}
\newcommand{\prob}{\mathbb P}

\newcommand{\giv}{\,|\,}

\newcommand{\eq}{\begin{equation}}
\newcommand{\en}{\end{equation}}
\newcommand{\FM}{{\rm FM}}

\newcommand{\re}[1]{\mbox{(\ref{#1})}}
\newcommand{\rem}[1]{\mbox{\rm (\ref{#1})}}

\newcommand{\Nat}{\Bbb N}

\def\endpf{\hfill $\Box$ \vskip0.5cm}
\def \proof{\noindent{\em Proof.\ }}

\theoremstyle{plain}

\newtheorem{theorem}{\large Theorem}
\newtheorem{proposition}[theorem] {\large Proposition}
\newtheorem{definition}[theorem]{\large Definition}

\newtheorem{lemma}[theorem]{\large Lemma}

\begin{document}

\title{Exchangeable partitions derived from Markovian coalescents with simultaneous multiple collisions
\thanks{Research supported by N.S.F. Grant DMS-0405779}}
\author{Rui Dong
\thanks{Department of Statistics, University of California, Berkeley; e-mail: ruidong@stat.Berkeley.EDU}  }

\date{June 11, 2007}
\maketitle

\begin{abstract}
\noindent

Kingman derived the Ewens sampling formula for random partitions 
from the genealogy model defined by a Poisson process of mutations 
along lines of descent governed by a simple coalescent process. 
M\"ohle described the recursion which determines 
the generalization of the Ewens sampling formula when
the lines of descent are governed by a coalescent with multiple collisions. 
In \cite{donggnedinpitman06} authors exploit an analogy with the 
theory of regenerative composition and partition structures, 
and provide various characterizations of the associated exchangeable 
random partitions. 
This paper gives parallel results for the further generalized model with 
lines of descent following a coalescent with simultaneous multiple collisions. 
\end{abstract}



\section{Introduction}

Given a large population with many generations, 
we track backward in time the family history of each individual 
in the current generation. 
As we track further, the family lines coalesce with each other, 
eventually all terminating at a common ancestor of current generation. 
This model used on biological study for genealogy of haploid type \cite{canning74} 
is the prototype of Kingman's theory of random coalescent processes \cite{kingman82a,  kingman82c, kingman82b}. 
In Kingman's coalescent \cite{kingman82a}, each collision only involves two parts. 
The idea is extended to coalescent with multiple collisions in \cite{lambda, sagitov} 
where every collision can involve two or more parts.  
This model is further developed into the theory of coalescent with 
simultaneous multiple collisions in \cite{moehlesagitove01, jason00}. 
See \cite{bertoingold04, boltszni98, donggoldmartin05, evanspitman98, sagitov03, tavare84, watter84} 
for related developments.

\par
Kingman \cite{kingman82b} indicated a connection between 
random partitions and coalescent processes. 
Suppose in above haploid model the family line of current generation 
is modeled by Kingman's coalescent, 
and the mutations are applied along the family lines according to 
a Poisson process with rate $\theta/2$ for some non-negative real number $\theta$. 
Define a partition by saying that two individuals are in the 
same block if there is no mutation along their family lines before they coalesce. 
Then the resulting random partition is governed by the Ewens sampling formula 
with parameter $\theta$. 
See  \cite[Section 5.1, Exercise 2]{CSP} and \cite{BertoinBook, Nordborg} 
for review and more on this idea. 
Recently, M{\"o}hle \cite{moehle1} applied this idea to 
the genealogy tree modeled by coalescents with multiple collisions 
and simultaneous multiple collisions. 
He studied the resulting family of partitions, 
and derived a recursion which determines them. 
See \cite{moehle2} for more properties of this family.

Dong, Gnedin and Pitman \cite{donggnedinpitman06} offered a different approach 
to the family of random partitions generated by Poisson marking along 
the lines of descent modeled by a coalescent with multiple collisions. 
In their work, 
each part of partitions is assigned one of two possible states: active or frozen, 
and a new class of continuous time partition-valued coalescent processes, 
namely {\em coalescents with freeze}, is introduced. 
Every coalescent with freeze has a terminal state with all blocks frozen, 
called the {\em final partition} of this process, 
whose distribution is characterized by M{\"o}hle's recursion \cite{moehle1}. 
In the spirit of \cite{RCS, RPS}, 
the authors studied the discrete time chains embedded in the coalescent with freeze, 
and from the consistency of their transition operators 
they derived a backward recursion satisfied by the decrement matrix, 
analogous to \cite[Theorem 3.3]{RCS}. 
This decrement matrix determines the final partition through M\"ohle's recursion. 
An integral representation for the decrement matrix was derived 
using algebraic methods as in \cite{RCS}. 
Moreover, adapting an idea from \cite{RPS}, 
the authors established a uniqueness result by constructing another Markov chain, 
with state space the set of partitions of a finite set, 
whose unique stationary distribution is the law of the final partition 
restricted to this set.

\par
As noted in \cite{donggnedinpitman06}, 
the current paper serves as a supplement with the theory on more general case provided. 
We focus on the family of partitions generated by Poisson marking along 
the lines of descent modeled by a Markov coalescent with simultaneous multiple collisions. 
This family was characterized by the generalized form of M\"ohle's recursion \cite{moehle1}, 
while in this paper, parallel analysis as that in \cite{donggnedinpitman06} is exploited.

\par
Notations and background are introduced in Section 2, 
together with a review of M\"ohle's idea. 
In Section \ref{defxicoalfreeze}, 
the generalized coalescent with freeze is defined, 
followed by the connection between this process and the generalized form of M\"ohle's recursion. 
Then in Section \ref{generalFM} 
we study the generalized freeze-and-merge (FM) operators 
of the embedded finite discrete chain of the coalescent with freeze process, 
those consistency with sampling gives a backward recursion for the 
generalized decrement matrix. 
Also our main result regarding finite partitions is stated. 
In Section \ref{generalSA} 
another partition-valued Markov chain with generalized sample-and-add (SA) operation is introduced, 
the law of the partition in our study is identified as the unique stationary distribution of this chain. 
Finally in Section 6 
we derive the integral representation for the generalized infinite decrement matrix, 
our main result regarding infinite partitions is stated. 
Section 3-6 of this paper can be seen as generalizations of 
the corresponding sections of \cite{donggnedinpitman06}.

\section{Notations and background}

Following the notations of \cite{donggnedinpitman06} and \cite{CSP}, 
a partition of a finite set $F$ into $\ell$ blocks, also called a {\em finite set partition}, 
is an unordered collection of non-empty disjoint sets $\{A_1,\ldots,A_\ell\}$ whose union is $F$. 
Partitions of the set $[n]:=\{1,2,\ldots,n\}$ for $n\in \mathbb{N}$ are of our special interest. 
Let $\mathcal{P}_{[n]}$ be the set of all partitions of $[n]$.
For a positive integer $n$, 
a {\em composition} of $n$ is an ordered sequence of positive integers 
$(n_1,n_2,\ldots,n_\ell)$ with $\sum_{i=1}^\ell n_i=n$, 
where $\ell\in \mathbb{N}$ is number of parts. 
Let $\mathcal{C}_n$ be the set of all compositions of $n$, and 
we use $\mathcal{P}_n$ to denote the set of non-increasing compositions of $n$, 
also called {\em partitions of $n$}.

\par
Take $\pi_n=\{A_1,A_2,\ldots,A_\ell\}$ as a generic partition of $[n]$, 
which we write as $\pi_n\vdash [n]$.
The {\em shape function} from partitions of the set $[n]$ to 
partitions of the positive integer $n$ is defined by
\eq
\mathtt{shape}(\pi_n)=(|A_1|,|A_2|,\ldots,|A_\ell|)^{\downarrow}
\en 
where $|A_i|$ represents the size of block $A_i$ which is the number of elements in the block,
 and ``$\downarrow$'' means arranging the sequence of sizes in non-increasing order.

 \par
A random partition $\Pi_n$ of $[n]$ is a random variable taking values in $\mathcal{P}_{[n]}$. 
It is called exchangeable if its distribution is invariant under the action 
on partitions of $[n]$ by the symmetric group of permutations of $[n]$. 
Equivalently, the distribution of $\Pi_n$ is then given by the formula
\eq
\mathbb{P}(\Pi_n=\{A_1,A_2,\ldots,A_\ell\}) = p_n(|A_1|, |A_2|,\ldots,|A_\ell|)
\en
for some symmetric function $p_n$ of compositions of $n$. 
$p_n$ is called the {\em exchangeable partition probability function (EPPF)} of $\Pi_n$.

\par
An {\em exchangeable random partition of $\Nat$} is a sequence of 
exchangeable set partitions $\Pi_{\infty}=(\Pi_n)_{n=1}^\infty$ with $\Pi_{n}\vdash [n]$, 
subject to the {\em consistency condition}
\eq\label{infinitepartition}
\Pi_n |_{m} = \Pi_m,
\en
where the restriction operator $|_{m}$ acts on $\mathcal{P}_{[n]}$, $n>m$, 
by deleting elements $m+1,m+2,\ldots,n$. 
The distribution of such an exchangeable random partition of $\Nat$ is determined by 
the function $p$ defined on the set of all integer compositions 
$\mathcal{C}_\infty:=\cup_{i=1}^\infty \mathcal{C}_i$, 
which coincides with the EPPF $p_n$ of $\Pi_n$ when acting on $\mathcal{C}_n$. 
This function $p$ is called the infinite EPPF associated with $\Pi_{\infty}$.
The consistency condition \re{infinitepartition} translates into 
the following {\em addition rule} for the EPPF $p$: 
for each positive integer $n$ and each composition $(n_1,n_2,\ldots,n_\ell)$ of $n$,
\eq\label{eppfaddition}
p(n_1,n_2, \ldots,n_\ell) = p(n_1,n_2, \ldots, n_\ell, 1)+\sum_{i=1}^\ell p(n_1, \ldots, n_i+1, \ldots, n_\ell)
\en
where $(n_1, \ldots, n_i+1, \ldots, n_\ell)$ is formed from $(n_1, \ldots, n_\ell)$ by adding $1$ to $n_i$. 
Conversely, if a nonnegative function $p$ on compositions satisfies (\ref{eppfaddition}) and 
the normalization condition $p(1)=1$,
then by Kolmogorov's extension theorem there exists an exchangeable random partition $\Pi_\infty$ with EPPF $p$.

\par
Similar definitions apply to a finite sequence of consistent exchangeable random set partitions 
$(\Pi_m)_{m=1}^n$ with $\Pi_{m}\vdash [m]$, where $n$ is some fixed positive integer. 
The {\em finite} EPPF $p$ of such a sequence can be defined as the unique recursive extension of $p_n$ 
by the addition rule \re{eppfaddition} to all compositions $(n_1,n_2,\ldots,n_\ell)$ of $m < n$.

\par
Let $\mathcal{P}_{\infty}$ be the set of all partitions of $\Nat$. 
We identify each $\pi_\infty \in \mathcal{P}_{\infty}$ as the sequence 
$(\pi_{1},\pi_{2},\ldots)\in \mathcal{P}_{[1]}\times \mathcal{P}_{[2]}\times\cdots$, 
where $\pi_{n}=\pi_\infty|_n$ is the restriction of $\pi_{\infty}$ to $[n]$ by 
deleting all elements bigger than $n$. 
Give $\mathcal{P}_{\infty}$ the topology it inherits as a subset 
of $\mathcal{P}_{[1]}\times \mathcal{P}_{[2]} \times\cdots$ with the product of discrete topologies, 
so the space $\mathcal{P}_{\infty}$ is compact and metrizable. 
Following \cite{evanspitman98, kingman82a,lambda}, 
call a $\mathcal{P}_{\infty}$-valued stochastic process $(\Pi_{\infty}(t), t\ge 0)$ a {\em coalescent} if 
it has c\`adl\`ag paths and $\Pi_{\infty}(s)$ is a {\em refinement} of $\Pi_{\infty}(t)$ for every $s<t$. 
For a non-negative finite measure $\Lambda$ on the Borel subsets of $[0,1]$, 
a {\em $\Lambda$-coalescent} is a $\mathcal{P}_{\infty}$-valued Markov coalescent 
$(\Pi_{\infty}(t), t\ge 0)$ whose restriction $(\Pi_{n}(t), t\ge 0)$ to $[n]$ is for each $n$ 
a Markov chain such that when $\Pi_{n}(t)$ has $b$ blocks, 
each $k$-tuple of blocks of $\Pi_{n}(t)$ is merging to form a single block at rate $\lambda_{b,k}$, where
\eq\label{lambdaint}
\lambda_{b,k}=\int_0^1x^{k-2}(1-x)^{b-k}\Lambda(dx)\ \ \ \ \ \ (2\le k\le b<\infty).
\en
When $\Lambda=\delta_0$, this reduces to Kingman's coalescent \cite{kingman82a, kingman82b, kingman82c} 
with only binary merges. 
When $\Lambda$ is the uniform distribution on $[0,1]$, 
the coalescent is the Bolthausen-Sznitman coalescent \cite{boltszni98}.

\par
A key property of the $\Lambda$-coalescent is that 
the collision rates do not depend on the internal structure of each block, 
it is therefor natural to seek a more general class of coalescent processes which 
retain this property and undergo ``silmutaneous multiple collisions''. 
This idea is mentioned in \cite[Section 3.3]{lambda}, 
and the coalescent process with simultaneous multiple collisions is 
first obtained and characterized in \cite{moehlesagitove01} by a sequence of measures, 
then in \cite{jason00} Schweinsberg finds a more compact characterization 
by a single non-negative measure $\Xi$ on the infinite simplex 
\eq
\Delta = \{(x_1,x_2,\ldots): x_1\ge x_2\ge \cdots\ge 0, \sum_{i=1}^\infty x_i\le 1\}.
\en
See \cite{sagitov03, CSP, moehle1} for various discussions for this process.

\par
To be more specific, following notations in \cite{jason00} 
for generic set partition with $b\in \Nat$ blocks, 
let $(k_1,k_2,\ldots, k_r, s)$ be a sequence of positive integers with 
$s\ge 0$, $r\ge 1$, $k_i\ge 2$ for $i=1,2,\ldots,r$, and $s+\sum_{i=1}^r k_i=b$, 
we define a {\em $(b; k_1,k_2,\ldots,k_r; s)$-collision} to be a merge of 
$b$ blocks into $r+s$ blocks in which $s$ blocks remain unchanged and 
the other $r$ blocks contain $k_1,k_2,\ldots,k_r$ of the original blocks. 
The order of $k_1,k_2,\ldots,k_r$ does not matter. 
For example, take the original partition as $\{ \{1,3\}, \{2\}, \{4\}, \{5\},\{6,7\},\{8\}\}$, 
then partition $\{\{1,2,3,5\},\{4,6,7\}, \{8\}\}$ is a $(6;3,2;1)$-collision of the original partition 
and also a $(6;2,3;1)$-collision of the original one. 
It is clear that for any generic set partition with $b$ blocks, 
the number of possible $(b; k_1,k_2,\ldots,k_r; s)$-collisions is
\eq\label{dnumber}
d(b;k_1,k_2,\ldots, k_r; s):=\frac{b!}{s!\ \Pi_{j=2}^b(j!)^{l_j}l_j!}= {b \choose {k_1\ \cdots\ k_r\ s}}\frac{1}{\Pi_{j=2}^b l_j!}
\en
where $l_j:=\#\{i: k_i=j\}$.

\par
Let $\Xi$ be some non-negative finite measure on the infinite simplex 
with the form $\Xi=\Xi_0+a\delta_0$, 
where $\Xi_0$ has no atom at zero and $\delta_0$ is a unit mass at zero. 
A {\em $\Xi$-coalescent} starting from generic infinite partition $\pi_\infty\in \mathcal{P}_\infty$ 
is a $\mathcal{P}_\infty$-valued coalescent $(\Pi_{\infty}(t), t\ge 0)$ with 
\begin{itemize}
\item $\Pi_\infty(0)=\pi_\infty$,
\item for each positive integer $n$, the restricted process 
$(\Pi_{n}(t), t\ge 0):=(\Pi_{\infty}(t), t\ge 0)|_n$ is a $\mathcal{P}_{[n]}$-valued 
Markov chain such that when $\Pi_n(t)$ has $b$ block, 
each possible $(b;k_1,k_2,\ldots, k_r; s)$-collision happens with rate $\lambda_{b; k_1,k_2,\ldots,k_r; s}$, 
which is defined as the integral
\eq\label{xilambdaint}
\int_{\Delta}\left . \left (\sum_{l=0}^s\sum_{i_1\ne \cdots \ne i_{r+l}}{s\choose l}x_{i_1}^{k_1}\cdots x_{i_r}^{k_r} x_{i_{r+1}}\cdots x_{i_{r+l}}(1-\sum_{j=1}^\infty x_j)^{s-l}\right)\right  / \sum_{j=1}^\infty x_j^2 \ \Xi_0(dx)+a 1_{\{r=1,k_1=2\}}
\en
\end{itemize}
It is called a standard $\Xi$-coalescent if the starting state $\pi_\infty=\Sigma_\infty$ 
being the partition of $\Nat$ into singletons.

\par
The measure $\Xi$ which characterizes the coalescent is derived from the consistency requirement, 
that is for any positive integers $0<m<n<\infty$, and $\pi_n\vdash [n]$, 
the restricted process $(\Pi_n(t)|_{m},\,t\ge 0)$ given $\Pi_{n}(0)=\pi_{n}$ 
has the same law as $(\Pi_{m}(t), t\ge 0)$ given $\Pi_{m}(0)=\pi_{n}|_{m}$. 
This condition is fulfilled if and only if the array of rates $(\lambda_{n; k_1,k_2,\ldots,k_r; s})$ satisfies
\eq\label{xilambdaconsis}
\lambda_{n; k_1,k_2,\ldots,k_r; s}=\sum_{i=1}^r \lambda_{n; k_1,\ldots,k_{i-1},k_i+1,k_{i+1},\ldots,k_r; s}
+s\lambda_{n+1; k_1,k_2,\ldots,k_r,2; s-1}+\lambda_{n+1; k_1,k_2,\ldots,k_r; s+1}
\en
where we say $s\lambda_{n+1; k_1,k_2,\ldots,k_r,2; s-1}=0$ when $s=0$ although it is undefined, 
so that the right hand side makes sense. 
The integral representation (\ref{xilambdaint}) can be derived from \re{xilambdaconsis} and 
exchangeability arguments \cite{jason00}.

\par
M{\"o}hle \cite{moehle1} studied the following generalization of Kingman's model \cite{kingman82b}. 
Take a genetic sample of $n$ individuals from a large population and 
label them as $\{1,2,\ldots,n\}$. 
Suppose the ancestral lines of these $n$ individuals evolve by the rules of 
a $\Lambda$-coalescent ($\Xi$-coalescent in general), 
and that given the genealogical tree, 
whose branches are the ancestral lines of these individuals, 
mutations occur along the ancestral lines according to a 
Poisson point process with rate $\rho>0$.
The infinite-many-alleles model is assumed, 
which means that when a gene mutates,  a brand new type appears.
Define a random partition of $[n]$ by declaring individuals $i$ and $j$ 
to be in the same block if and only if they are of the same type,
that is either $i=j$ or there are no mutations along the ancestral lines 
of $i$ and $j$ before these lines coalesce. 
These random partitions are exchangeable, and consistent as $n$ varies. 
The EPPF of this random partition is the unique solution $p$ with $p(1) = 1$ of 
{\em \Mohle's recursion}. 
In this paper, we focus on the general case with the ancestral lines modeled by a $\Xi$-coalescent.

\par
In order to write out M\"ohle's recursion \cite[Theorem 5.1]{moehle1} of the general case 
in a form which fits our treatment better, we introduce some notations: 
given a generic composition $(n_1,n_2,\ldots,n_\ell)$ of positive integer $n$ and 
a positive integer set $\{k_1,\ldots,k_r\}$ with $r\ge 1, k_1,\ldots,k_r\ge 2, \sum_{j=1}^r k_j\le n$, 
let each $n_i$ choose some $k_j$'s from $\{k_1,\ldots,k_r\}$. 
Denote the index set of those $k_j$'s chosen by $n_i$ as 
$$
\eta_i:=\{j:\ k_j\ \mbox{is chosen by}\ n_i\},
$$
and it is $\O$ when $n_i$ chooses nothing. 
The choices must satisfy
\begin{itemize} 
\item every $k_j$ can only be chosen by one $n_i$; 
\item for each $i$, $n_i\ge \sum_{j\in \eta_i}k_j$. 
\end{itemize}
So every such choice can be represented by a sequence of index sets 
\eq
\eta=(\eta_1,\ldots,\eta_\ell),
\en
We write elements in $\eta_i$ as $\eta_{i(1)},\eta_{i(2)},\ldots$, 
and use the notation
\eq
H_{\{k_1,\ldots,k_r\}}^{(n_1,\ldots,n_\ell)}=\{\eta^1,\eta^2,\ldots\}
\en
to denote the set of all possible choices, 
especially we identify two choices $\eta^1$ and $\eta^2$ if for each $i$, 
$\{k_j: j\in \eta^1_i\}=\{k_j: j\in \eta^2_i\}$. 
This eliminates some trivial repetition caused by the indexing of  $\{k_1,\ldots,k_r\}$. 
For example, in our definition $H_{\{3,3,3\}}^{(6,3)}$ contains only one choice $(\{1,2\},\{3\})$, 
which is considered as the same as $(\{2,3\},\{1\})$ or $(\{1,3\},\{2\})$.

With these notations, we write M\"ohle's recursion as
\begin{align}\label{generalrec1}
&p(n_1,n_2,\ldots,n_\ell)=
\frac{q(n:1)}{n}\sum_{j:n_j=1}\,p(\widehat{n_j})+
\sum_{\{k_1,\ldots,k_r\}}q(n:k_1,\ldots,k_r;n-\sum_{j=1}^rk_j)\,\nonumber\\
& \sum_{\eta\in H_{\{k_1,\ldots,k_r\}}^{(n_1,\ldots,n_\ell)}} 
\frac{\prod_{i=1}^\ell d(n_i; k_{\eta_{i(1)}}, \ldots, k_{\eta_{i(| \eta_i |)}};  n_i-\sum_{l=1}^{|\eta_i|} k_{\eta_{i(l)}} ) }
{d(n; k_1,\ldots,k_r;n-\sum_{j=1}^rk_j)} 
p(n_1-\sum_{l=1}^{|\eta_1|} k_{\eta_{1(l)}} + |\eta_1|,\ldots, n_\ell-\sum_{l=1}^{|\eta_\ell|} k_{\eta_{\ell(l)}}+|\eta_\ell| )
\end{align}
where $(\widehat{n_j})$ is formed from $(n_1,n_2,\ldots,n_\ell)$ by deleting part $n_j$, 
$|\eta_i|$ is the number of elements in $\eta_i$ and 
the second sum on the right hand side is over all 
integer sets $\{k_1,\ldots,k_r\}$ with $r\ge 1, k_1,\ldots,k_r\ge 2, \sum_{i=1}^r k_i\le n$. 
Also in \re{generalrec1}, 
\eq\label{xiqsequencek}
q(b; k_1,k_2,\ldots,k_r; s):=\frac{\Phi(b; k_1,k_2,\ldots,k_r; s)}{\Phi(b)},
\en
\eq\label{xiqsequence1}
 q(b:1):={\Phi(b:1)\over \Phi(b)}.
\en
where
\begin{eqnarray}
\Phi(b:1)&:=&\rho b\\
\Phi(b; k_1,k_2,\ldots,k_r; s)&:=&d(b;k_1,k_2,\ldots, k_r; s)\lambda_{b; k_1,k_2,\ldots,k_r; s}\\
\Phi(b)&:=&\Phi(b:1)+\sum_{\{k_1,\ldots,k_r\}}\Phi(b; k_1,k_2,\ldots,k_r; s)
\label{Phi33}.
\end{eqnarray}
in which the last sum is over all multisets $\{k_1,\ldots,k_r\}$ with 
$r\ge 1, k_1,\ldots,k_r\ge 2, \sum_{i=1}^r k_i\le b$, 
and $\lambda_{b; k_1,k_2,\ldots,k_r; s}$'s are defined as the integral \re{xilambdaint}. 
The meaning of these notations are obvious: 
suppose at some time $t\ge 0$, 
a finite $\Xi$-coalescent freezing at rate $\rho$ has $b$ active blocks, 
then $\Phi(b:1)$ is the total rate at which some active block freezes; 
$\Phi(b; k_1,k_2,\ldots,k_r; s)$ is the total rate of a $(b; k_1,k_2,\ldots,k_r; s)$-collision, 
and $\Phi(b)$ is the total rate of changing state. 
If we look at its embedded discrete jump chain, 
$q(b:1)$, $q(b; k_1,k_2,\ldots,k_r; s)$ are the transition probabilities 
of these two kinds of events, respectively.

\par
M{\"o}hle \cite{moehle1} derived the recursion \re{generalrec1} 
by conditioning on whether the first event met tracing back in time 
from the current generation is a mutation or collision. 
On the left side of \re{generalrec1},
$p(n_1,n_2,\ldots,n_\ell)$ is the probability of ending up with any {\em particular} partition 
$\pi_n$ of the set $[n]$ into $\ell$ blocks of sizes $(n_1,n_2,\ldots,n_\ell)$.
On the right side, $q(n:1)$ is the chance that starting from the current generation, 
one of the $n$ genes mutates before any collision; 
for this to happen together with the specified partition of $[n]$, 
the individual with this gene must be chosen from those among the singletons of $\pi_n$, 
with chance  $1/n$ for each different choice, 
and after that the restriction of the coalescent process to a subset of $[n]$ of 
size $n-1$ must end up generating the restriction of $\pi_n$ to that set.
Similarly, $q(n:k_1,\ldots,k_r;n-\sum_{j=1}^rk_j)$ is the chance that the first event met is 
a $(n:k_1,\ldots,k_r;n-\sum_{j=1}^rk_j)$-collision. 
By the definition, 
each element $\eta$ in $H_{\{k_1,\ldots,k_r\}}^{(n_1,\ldots,n_\ell)}$ represents 
a class of ways of grouping singleton blocks to perform a $(n:k_1,\ldots,k_r;n-\sum_{j=1}^rk_j)$-collision 
such that it is possible for the resulting partition to have block sizes $\{n_1,n_2,\ldots,n_\ell\}$. 
Given each $\eta$, there are still various possibilities due to the different grouping scheme inside 
each block with sizes $n_i$, this is where the factor 
\eq
\frac{\prod_{i=1}^\ell d(n_i; k_{\eta_{i(1)}}, \ldots, k_{\eta_{i(| \eta_i |)}};  n_i-\sum_{l=1}^{|\eta_i|} k_{\eta_{i(l)}} ) }
{d(n; k_1,\ldots,k_r;n-\sum_{j=1}^rk_j)} 
\en
in \re{generalrec1} comes from. 
Conditioning on a particular selection, 
the restriction of the coalescent process to some set of $n-\sum_{i=1}^rk_i+r$ lines of descent 
must end up generating a particular partition of these $n-\sum_{i=1}^rk_i+r$ lines into sets of sizes 
\eq
(n_1-\sum_{l=1}^{|\eta_1|} k_{\eta_{1(l)}} + |\eta_1|,\ldots, n_\ell-\sum_{l=1}^{|\eta_\ell|} k_{\eta_{\ell(l)}}+|\eta_\ell| ).
\en
The multiplication of various probabilities is justified by the strong Markov property 
of the $\Xi$-coalescent at the time of the first event,
and by the special symmetry property that lines of descent representing blocks of individuals 
coalesce according to the same dynamics as if they were singletons.

\par 
Same as \cite{donggnedinpitman06}, 
in this paper we choose to step back from the special forms \re{xiqsequencek}, \re{xiqsequence1} 
of the $q$ entries $q(b; k_1,k_2,\ldots,k_r; s)$'s and $q(n:1)$ 
derived from the $(\Xi, \rho)$, 
and analyse M{\"o}hle's recursion \re{generalrec1} as an abstract relation between 
an array of $q$ entries and a function of compositions $p$. 
In particular, we ask the following questions, quoted from \cite{donggnedinpitman06}:
\begin{enumerate}
\item For which array of $q$ entries is \Mohle's recursion \re{generalrec1} 
satisfied by the EPPF $p$ of some exchangeable random partition of $[n]$, 
and is this $p$ uniquely determined?
\item How can such random partitions be characterized probabilistically?
\item Can such random partitions of $[n]$ be consistent as $n$ varies for 
any other array of $q$ entries besides those derived from $(\Xi, \rho)$ as above?
\end{enumerate}
We stress that in the first two questions 
the recursion \re{generalrec1} is only required to hold for a single value of $n$, 
while in the third question \re{generalrec1} must hold for all $n=1,2,\ldots$.
The answer to the first question is that 
for each fixed array of $q$ entries with sum equaling to $1$, 
which will be made precise later in Section \ref{generalFM},  
\Mohle's recursion \re{generalrec1} determines a unique EPPF $p$ 
for an exchangeable random partition of $[n]$ (Theorem \ref{finitegeneralrec}). 
Answering the second question, 
we characterize the distribution of this random partition in two different ways: 
firstly as the terminal state of a discrete-time Markovian coalescent process, 
the generalized {\em freeze-and-merge chain} introduced in Section \ref{generalFM}, 
and secondly as the stationary distribution of a partition-valued Markov chain 
with quite a different transition mechanism, 
the generalized {\em sample-and-add chain} introduced in Section \ref{generalSA}. 
The answer to the third question is positive if we restrict $n$ to some bounded range of values, 
for some but not all $q$ (see Section \ref{generalFM}), 
but negative if we require consistency for all $n$ (Theorem \ref{generalmain}): 
if an infinite EPPF $p$ solves \Mohle's recursion \re{generalrec1} for all
$n$ for some non-negative $q$ entries, 
then $q$ entries must have the form \re{xiqsequencek}, \re{xiqsequence1} for some $(\Xi,\rho)$. 
All these results are generalized version of the corresponding ones in \cite{donggnedinpitman06}.

\par
The analysis in this paper follows the same route as \cite{donggnedinpitman06}, 
where the authors were guided by a remarkable parallel between 
the theory of finite and infinite partitions subject to \Mohle's recursion and 
the theory of {\em regenerative partitions} developed in \cite{RCS, RPS}. 
Many of these parallels are summarized in \cite[Section 9]{donggnedinpitman06}.
Readers can check \cite{donggnedinpitman06} for other aspects of this idea.

\section{Coalescents with freeze}\label{defxicoalfreeze}

Same as \cite{donggnedinpitman06}, 
we consider the structure of a partition of a set (respectively, of an integer) with 
each of its blocks (or parts) assigned one of two possible conditions, 
which we call {\em active} and {\em frozen}.  
We name such a combinatorial object a {\em partially frozen partition}  of a set or of an integer. 
Use symbol $\Sigma_n^*$ for the pure singleton partition of $[n]$ 
with all blocks active, 
and $\Sigma_\infty^*$ for the sequence $(\Sigma_n^*)_{n=1}^\infty$. 
We include the possibilty of all blocks being active or frozen as 
special cases of partially frozen partitions.
Ignoring the conditions of the blocks of a partially frozen partition  $\pi^*$ 
{\em induces} an ordinary partition $\pi$. 
The {\em *-shape} of a partially frozen partition $\pi_n^*$ of $[n]$ is 
the corresponding partially frozen partition of $n$, 
and the ordinary shape is defined in terms of the induced partition $\pi_n$.

\par
For each positive integer $n$, 
we denote  $\mathcal{P}_{[n]}^*$ the set of all partially frozen partitions of $[n]$. 
Let $\mathcal{P}_{\infty}^*$ be the set of all partially frozen partitions of $\Nat$. 
We identify each element $\pi_\infty^*\in \mathcal{P}_{\infty}^*$ as the sequence 
$(\pi_1^*,\pi_2^*,\ldots)\in \mathcal{P}_{[1]}^*\times \mathcal{P}_{[2]}^*\times\cdots$, 
where $\pi_n^*$ is $\pi_\infty^*|_{n}$ the restriction of $\pi_\infty^*$ to $[n]$. 
Endowing $\mathcal{P}_{\infty}^*$ with the topology it inherits as a subset of 
$\mathcal{P}_{[1]}^*\times \mathcal{P}_{[2]}^*\times\cdots$, 
the space $\mathcal{P}_{\infty}^*$ is compact and metrizable. 
We call a random partially frozen partition of $[n]$ exchangeable
if its distribution is invariant under the action of permutations of $[n]$. 
Following \cite{donggnedinpitman06}, 
call a $\mathcal{P}_{\infty}^*$-valued stochastic process 
$(\Pi_{\infty}^*(t), t\ge 0)$ a {\em coalescent} if it has c\`adl\`ag paths and 
$\Pi_{\infty}^*(s)$ is a {\em *-refinement} of $\Pi_{\infty}^*(t)$ for every $s<t$, 
meaning that the induced partition $\Pi_{\infty}(s)$ is a refinement of $\Pi_{\infty}(t)$ 
and the set of frozen blocks of $\Pi_{\infty}^*(s)$ is a subset of 
the set of frozen blocks of $\Pi_{\infty}^*(t)$.

\par
The remaining part of the section gives formal statement about the  
{\em $\Xi$-coalescent with freeze}, 
which is the generalization of $\Lambda$-coalescent with freeze 
defined in \cite{donggnedinpitman06}. 
The connection between M\"ohle's model \cite{moehle1} and our realization 
has been outlined clearly in \cite[Section 3]{donggnedinpitman06}.

\begin{theorem}\label{Xicoalescentwithfreeze}
Let $\{ \lambda_{b; k_1,k_2,\ldots,k_r; s}: 2\le b<\infty, r\ge 1, k_1,\ldots,k_r\ge 2,s\ge 0,b=s+\sum_{i=1}^r k_i\}$, 
$(\rho_{n,b}, 1\le b\le n<\infty)$ be two arrays of non-negative real numbers. 
There exists for each $\pi_\infty^*\in \mathcal{P}_{\infty}^*$ 
a $\mathcal{P}_{\infty}^*$-valued coalescent $(\Pi_{\infty}^*(t), t\ge 0)$ with 
$\Pi_{\infty}^*(0)=\pi^*_\infty$, 
for each $n$ whose restriction $(\Pi_{n}^*(t), t\ge 0)$ to $[n]$ is 
a $\mathcal{P}_{[n]}^*$-valued Markov chain starting from $\pi_n^*=\pi_\infty^*|_{n}$, 
and evolving with the rules:
\begin{itemize}
\item at each time $t \ge 0$, conditionally given $\Pi_{n}^*(t)$ with $b$ active blocks, 
each possible $(b; k_1,k_2,\ldots,k_r; s)$-collision is occurring with rate 
$\lambda_{b; k_1,k_2,\ldots,k_r; s}$, and 

\item each active block turns into a frozen block at rate $\rho_{n,b}$,
\end{itemize}
if and only if the integral representation \rem{xilambdaint} holds 
for some non-negative finite measure on the infinite simplex with the form $\Xi=\Xi_0+a\delta_0$, 
where $\Xi_0$ has no atom at zero and $\delta_0$ is a unit mass at zero; 
and $\rho_{n,b}=\rho$ for some non-negative real number $\rho$. 
This $\mathcal{P}_{\infty}^*$-valued process $(\Pi_\infty^*(t), t\ge 0)$ directed by
 $(\Xi,\rho)$ is a strong Markov process. 
 
For $\rho=0$, this process reduces to the $\Xi$-coalescent, 
and for $\rho >0 $ the process is obtained by superposing Poisson marks
at rate $\rho$ on the merger-history tree of a $\Xi$-coalescent, 
and freezing the block containing $i$ at the time of the first mark along 
the line of descent of $i$ in the merger-history tree.
\end{theorem}
\proof
Consistency of the rate descriptions for different $n$ implies \re{xilambdaconsis}, 
from which we have existence of measure $\Xi$ and the integral representation \re{xilambdaint} 
by \cite[Lemma 18]{jason00} and \cite[Theorem 2]{jason00}. 
Equality of the $\rho_{n,b}$'s is obvious by consistency.
\endpf

\begin{definition}\label{defxiwithfreeze}
{\rm 
Call this $\mathcal{P}_{\infty}^*$-valued Markov process directed by 
a non-negative integer $\rho$ and a non-negative finite measure $\Xi$ on the infinite simplex 
the {\em $\Xi$-coalescent freezing at rate $\rho$},
or the {\em $(\Xi,\rho)$-coalescent} for short.
Call a $(\Xi,\rho)$-coalescent starting from state $\Sigma_\infty^*$ a 
{\em standard $\Xi$-coalescent freezing at rate $\rho$}, 
where $\Sigma_\infty^*$ is the pure singleton partition with all blocks active.}
\end{definition}

Consider the finite coalescent with freeze $(\Pi^*_n(t), t\ge 0)$ which is 
the restriction of a standard $\Xi$-coalescent freezing at rate $\rho$ to $[n]$. 
It is clear that as long as the freezing rate $\rho$ is positive, 
in finite time the process $(\Pi_n^*(t),t\ge 0)$ will eventually reach a 
{\em final partition} $E_n^*$, with all of its blocks in the  frozen condition. 
Set $E_\infty^*:=(E_n^*)$ as the final partition of $(\Pi_\infty^*(t), t\ge 0)$, 
and denote its induced partition as $E_\infty=(E_n)$. 
If we look at the discrete chain embedded in finite $\Xi$-coalescent freezing at rate $\rho$, 
by conditioning on the first transition we can see the following facts:

\begin{theorem} {\em ( M\"ohle \cite[Theorem 5.1]{moehle1})}
The induced final partition $E_\infty=(E_n)_{n=1}^\infty$ of 
a standard $\Xi$-coalescent freezing at rate $\rho>0$ is 
an exchangeable infinite random partition of $\Nat$ whose
EPPF $p$ is the unique solution of \Mohle's recursion {\rm \re{generalrec1}} with 
$q$ coefficients defined through $(\Xi,\rho)$ as in {\rm \re{xiqsequencek}, \re{xiqsequence1}}. 
\end{theorem}

\section{Freeze-and-merge operations.}\label{generalFM}

Following \cite[Section 4]{donggnedinpitman06}, 
given a continuous time stochastic process $X$ with right continuous piecewise constant path, 
the {\em jumping process derived from $X$} is the discrete-time process
$$
\widehat{X}=(\widehat{X}(0),\widehat{X}(1),\ldots) = ( X(T_0), X(T_1), X(T_2), \ldots )
$$
where $T_0 := 0$ and $T_k$ for $k \ge 1$ is the least $t > T_{k-1}$ such that $X(t) \neq X(T_{k-1})$,
if there is such a $t$, and $T_k = T_{k-1}$ otherwise. 
In particular, the finite coalescent with freeze $(\Pi^*_n(t), t\ge 0)$, obtained by
restriction to $[n]$ of a $\Xi$-coalescent freezing at positive rate $\rho$, 
is a Markov chain with transition rate 
$d(b;k_1,k_2,\ldots, k_r; s) \lambda_{b; k_1,k_2,\ldots,k_r; s}$ for a $(b; k_1,k_2,\ldots,k_r; s)$-collision 
and rate $b\rho$ for a freeze,
where $b$ is the number of active blocks at time $t$ and 
the $d(b;k_1,k_2,\ldots, k_r; s)$'s and $\lambda_{b,k}$'s are as in 
\re{dnumber} and \re{xilambdaint}; 
while the jumping process $\widehat{\Pi}^*_n$ is then a Markov chain 
governed by the following {\em freeze-and-merge} operation  ${\rm FM}_n$, 
which acts on a generic partially frozen partition $\pi^*_n$ of $[n]$ as follows:  
if $\pi^*_n$ has $b>1$ active blocks then
\begin{itemize}
\item with probability $q(b:k_1,k_2,\ldots, k_r; s)$ a selection of active blocks 
to perform a $(b; k_1,k_2,\ldots,k_r; s)$-collision is chosen uniformly at random 
from  $d(b;k_1,k_2,\ldots, k_r; s)$ total number of possible choices and 
the $(b; k_1,k_2,\ldots,k_r; s)$-collision is then performed as the chosen way;

\item with probability $q(b:1)$ an active block is chosen uniformly at random 
from $b$ blocks and turned into a frozen block.
\end{itemize}
where $q(b:k_1,k_2,\ldots, k_r; s)$'s and $q(b:1)$ are of the forms \re{xiqsequencek}, \re{xiqsequence1}. 
When $b=1$, only the second option is available. 
As a fact from last section, the continuous time processes $\Pi^*_n(t)$ 
are Markovian and consistent as $n$ varies, 
meaning that $\Pi^*_m(t)$ for $m <n$ coincides with 
$\Pi^*_n(t)|_m$, the restriction of $\Pi^*_n(t)$ to $[m]$.

\par 
To view \Mohle's recursion \re{generalrec1} in greater generality, 
we consider this freeze-and-merge operation  ${\rm FM}_n$ for $n$ 
some fixed positive integer, and an array
\eq\label{generalqsequence}
q^{(n)}:=(q(1:\cdot),q(2:\cdot),\ldots,q(n-1:\cdot),q(n:\cdot))
\en
where
\eq
q(1:\cdot)=\{q(1:1)\}=\{1\}
\en
and for $1<b\le n$,
\eq\label{generalqb}
q(b:\cdot):=\{q(b:1), q(b:k_1,k_2,\ldots, k_r; s):\  r\ge 1, k_i\ge 2\  {\rm for}\  i=1,2,\ldots,r, \ {\rm and}\ s=b-\sum_{i=1}^r k_i\ge 0 \}
\en
with all entries added up to $1$, where the order of index entries $k_i$ is neglected. 
And we always assume $q(b:\cdot)$ include all $q(b:k_1,k_2,\ldots, k_r; s)$'s for all 
possible indexes $\{k_1,k_2,\ldots,k_r\}$ by supplementing $0$ entries.

\par
Let $(\widehat{\Pi}_n^*(k), k=0,1,2,\ldots)$ 
be the Markov chain obtained by iterating 
${\rm FM}_n$ starting from $\widehat{\Pi}_n^*(0)=\Sigma_n^*$.  
The array $q^{(n)}$ can be see seen as generalization of decrement matrix $q_n$ 
in \cite[Section 4]{donggnedinpitman06}. 
For completeness, 
we include the propositions and lemmas listed in \cite[Section 4]{donggnedinpitman06}, 
most of them do not rely on the merging mechanism of freeze-and-merge operations. 
In particular, Lemma \ref{generalqrec} is the key result for this generalized case, 
and provides the basis for Theorem \ref{finitegeneralrec}, 
our main result regarding finite partitions.

\par
Observe that for $m=1,\ldots,n$ the first $m$ entries
\eq
(q(1:\cdot),q(2:\cdot),\ldots,q(n-1:\cdot),q(m:\cdot))
\en
of array $q^{(n)}$ comprise another array which itself defines a 
freeze-and-merge operation ${\rm FM}_m$ on partially frozen partitions of $[m]$.

\begin{proposition}\label{fmrec}
Given an array $q^{(n)}$ as in {\rm \re{generalqsequence}} and $(\widehat{\Pi}_n^*(k), k=0,1,2,\ldots)$ 
as a Markov chain governed by ${\rm FM}_n$ starting from $\Sigma_n^*$

{\em (i)} The $(\widehat{\Pi}_n^*(k))$ chain is strictly transient, it finally reaches 
a partially frozen partition $E_n^*$ of $[n]$ with all blocks frozen. 
Same thing holds for Markov chains governed by ${\rm FM}_m$, $1\le m\le n$ 
derived from $q^{(n)}$, note their final state as $E_m^*$'s, respectively.
Let $E_m$ be the induced partition of $[m]$ from $E_m^*$ for $1\le m\le n$.

{\em (ii)} Define $p$ as the function on  $\cup_{m=1}^n{\cal C}_m$ whose restriction to 
${\cal C}_m$ is the EPPF of $E_m$. 
Then $p$ satisfies \Mohle's recursion {\rm (\ref{generalrec1})} for  each composition $(n_1,n_2,\ldots,n_\ell)\in {\cal C}_n$.
\end{proposition}

\par
Part (ii) in the proposition follows by conditioning 
on the first transition of the $(\widehat{\Pi}_n^*(k))$ chain, 
similar with \cite[Lemma 4.1]{donggnedinpitman06}.

\par
In the general settings of Proposition \ref{fmrec}, the sequence of exchangeable 
final partitions $(E_m)_{m=1}^n$ need not be consistent with 
respect to restrictions. 
The question of what constraints on $q^{(n)}$ can guarantee the consistency of $(E_m)_{m=1}^n$ 
guided our reasoning. 

\begin{definition}\label{defqconsis}
For an array $q^{(n)}$ as in {\rm \re{generalqsequence}}  and $1 \le m < n$, 
call the transition operators ${\rm FM}_n$ and ${\rm FM}_m$  derived from $q^{(n)}$
{\em consistent} if whenever $\widehat{\Pi}_n^*$ is a Markov chain governed by ${\rm FM}_n$, 
the jump process derived from the restriction of $\widehat{\Pi}_n^*$ to $[m]$ is 
a Markov chain governed by ${\rm FM}_m$.
Call the decrement matrix $q_{n}$ consistent if this condition holds for every $1 \le m < n$.
\end{definition}

\par
It is clear from consistency of the continuous time $\Xi$-coalescent with freeze 
$(\Pi_n^*(t),~ t \ge 0 )$ introduced in the last section 
that for every $n$ the
corresponding array $q^{(n)}$ with forms \re{xiqsequencek} \re{xiqsequence1} is consistent.
Let $\FM_n(\pi_n^*)$ denote the random partition obtained by action of $\FM_n$ on an initial
partially frozen partition $\pi_n^*$ of $[n]$,

\begin{lemma}\label{factsqconsis} 
Given a particular array $q^{(n)}$ as in {\rm \re{generalqsequence}}:

{\em (i)} For fixed $1 \le m < n$ the transition operators ${\rm FM}_m$ and ${\rm FM}_n$ are consistent if and only if for
each  partially frozen partition $\pi_n^*$ of $[n]$, there is the equality in distribution
\eq\label{altconsis}
\FM_m(\pi_n^*|_{m})\ed\FM_n(\pi_n^*)||_{m}
\en
where on the left side $\pi_n^*|_{m}$ is the restriction of $\pi_n^*$ to $[m]$,
and on the right side the notation $||_{m}$ means the restriction to $[m]$ {\em conditional} 
on the event $\FM_n(\pi_n^*|_{m})\neq \pi_n^*|_{m}$ that $\FM_n$  freezes or merges at least one of the blocks of $\pi_n^*$ 
containing some element of $[m]$.

{\em (ii) } If ${\rm FM}_{m-1}$ and ${\rm FM}_m$ are consistent for every $1 <m \le n$, then 
so are ${\rm FM}_{m}$ and ${\rm FM}_n$ for every $1 <m \le n$; that is, $q^{(n)}$ is consistent.
\end{lemma}

\par
Following is the consistency results for array $q^{(n)}$, 
which is the generalized version of  \cite[Lemma 4.4]{donggnedinpitman06}. 
The omitted proof uses quite the same idea as that of \cite[Lemma 4.4]{donggnedinpitman06}, 
by looking at ${\rm FM}_{b+1}$ and ${\rm FM}_b$ applied to $\Sigma_{b+1}^*$ and $\Sigma_{b}^*$, 
respectively, and utilizing equation \re{altconsis}. 
Lemma \ref{factsqconsis} then links us between relation \re{altconsis} and 
consistency of $q^{(n)}$.

\begin{lemma}\label{generalqrec}
An array $q^{(n)}$ with form {\rm \re{generalqsequence}} is consistent if and only if 
it satisfies the backward recursion: 
\begin{align}\label{generalqrec1}
q(b: k_1,\ldots,k_r;s)=
&\sum_{i=1}^r \frac{(k_i+1)(l_{k_i+1}+1)}{(b+1) l_{k_i}}q(b+1:k_1,\ldots,k_{i-1},k_i+1,k_{i+1},\ldots,k_r;s)\nonumber\\
&+\frac{2(l_2+1)}{b+1}q(b+1:k_1,\ldots,k_r,2;s-1) 
+\frac{s+1}{b+1}q(b+1: k_1,\ldots,k_r;s+1)\nonumber\\
&+\frac{1}{b+1}q(b+1:1)q(b: k_1,\ldots,k_r;s) +\frac{2}{b+1}q(b+1: 2; b-1)q(b: k_1,\ldots,k_r;s)\nonumber\\
&\ \ \ \ \ \ \ \ \ \ \ \ \ \ \ \ \ \ \ \ \ \ \ \ \ \ \ \ \ \ \ \ \ \ \ \ \ \ \ \ \ \ \ \ \ \ \ \ \ \ \ \ \ \ \ \ \ \ \ \ \ \ \ ~~~~~(2 \le b <n),~~~~~
\end{align}

\begin{align}\label{generalqrec1b}
q(b:1)=\frac{b}{b+1}q(b+1:1)+\frac{1}{b+1}q(b+1:1)q(b:1)+\frac{2}{b+1}q(b+1:2;b-1)q(b:1)
~~~~~(1\le b <n),~~~~~
\end{align}
where $l_j:=\#\{i: k_i=j\}$, and when $s=0$, we say $q(b+1:k_1,\ldots,k_r,2;s-1)=0$ even though 
it is undefined, so that the right hand side of {\rm \re{generalqrec1}} makes sense.

Consequently, each array $q(n:\cdot)$ with form {\rm \re{generalqb}} determines 
a unique consistent $q^{(n)}$.
\end{lemma}

\par
With Proposition \ref{fmrec}, we have

\begin{lemma}\label{fmfinal}
For $1 \le m \le n$, let $E_m$ be the final partition of the $\FM_n$-chain starting in state $\Sigma_m^*$ 
as defined in Proposition {\rm \ref{fmrec}}.
If the array $q^{(n)}$ is consistent then the finite sequence of exchangeable 
random set partitions $(E_m)_{m=1}^n$ is consistent in the sense that
$$E_m\ed E_n|_m\,.$$
The finite EPPF $p$ of $(E_m)_{m=1}^n$ then satisfies \Mohle's recursion
{\rm (\ref{generalrec1})} for all compositions of $m\le n$ in the left hand side.
\end{lemma}

\par 
Here is our principal result regarding finite partitions satisfying (\ref{generalrec1}), 
which is parallel with \cite[Theorem 4.6]{donggnedinpitman06}.

\begin{theorem}\label{finitegeneralrec} 
For a positive integer $n>1$ and arbitrary array $q(n:\cdot)$ with form {\rm \re{generalqb}}
\begin{itemize}
\item[{\rm(i)}]
there exists a unique finite EPPF $p$ for a consistent sequence of random set partitions 
$(\Pi_m)_{m=1}^n$ which satisfies \Mohle's recursion {\rm (\ref{generalrec1})}
for all compositions of $n$, 
\item[{\rm(ii)}]
this finite EPPF $p$ satisfies \Mohle's recursion {\rm(\ref{generalrec1})} for all 
compositions of positive integers $m<n$ with coefficient arrays $q(m:\cdot)$ 
derived from $q(n:\cdot)$ by the recursion {\rm (\ref{generalqrec1}), (\ref{generalqrec1b})},
\item[{\rm(iii)}]
for each $1 \le m \le n$ the distribution of $\Pi_m$ determined by the 
restriction of this EPPF $p$ to compositions  of $m$ is that of
the final partition of the $\FM_m$ Markov chain with array $q^{(m)}$ defined  by {\rm (ii)}, 
starting from state $\Sigma_m^*$.
\end{itemize}
\end{theorem}

\proof 
Given $q(n: \cdot)$ , we can define a consistent array $q^{(n)}$ by 
the backward recursion (\ref{generalqrec1}),(\ref{generalqrec1b}).
Then use $q^{(n)}$ to build a sequence of Markov chains: for each $m$, 
the chain $(\widehat{\Pi}_m(k),k=0,1,2,\ldots)$ starts from $\Sigma_m^*$ and 
evolves according to $\FM_m$. 
By Lemma \ref{fmfinal}, the  sequence of induced final partitions $(E_m)_{m=1}^n$ 
of these chains has EPPF $p$ which satisfies recursion (\ref{generalrec1}). 
Hence the existence part of (i) follows.
The uniqueness in part (i) can be read from results in the next section.
The assertions (ii) and (iii) follow directly from this construction.
\endpf

\section{The sample-and-add operation.}\label{generalSA}

Following \cite[Section 5]{donggnedinpitman06}, 
we give \Mohle's recursion {\rm(\ref{generalrec1})} another interpretation 
as the system of equations for the invariant probability measure of 
a particular Markov transition mechanism on partitions of $[n]$, 
generalized sample-and-add operations, 
which takes the operations in \cite[Section 5]{donggnedinpitman06} as special forms. 
As a consequence, the uniqueness in Theorem \ref{finitegeneralrec} follows from 
the uniqueness of this invariant probability distribution.

\par
Fix some positive integer $n$ and a sequence
\eq\label{sageneralqb}
q(n:\cdot):=\{q(n:1), q(n:k_1,k_2,\ldots, k_r; s):\  r\ge 1, k_i\ge 2\  {\rm for}\  i=1,2,\ldots,r, \ {\rm and}\ s=n-\sum_{i=1}^r k_i\ge 0 \}
\en
with all entries added up to $1$, where the order of index entries $k_i$ is neglected. 
And we always assume $q(n:\cdot)$ include all $q(n:k_1,k_2,\ldots, k_r; s)$'s for all 
possible indexes $\{k_1,k_2,\ldots,k_r\}$ by supplementing $0$ entries. 
Let $K_n$ be a random element with its distribution according to this sequence $q(n:\cdot)$, 
i.e. $K_n$ equals to $1$ with probability $q(n:1)$, and equals to set $\{k_1,\ldots,k_r\}$ with 
probability $q(n:k_1,k_2,\ldots, k_r; s)$.

\par 
Consider the following {\em sample-and-add} random operation on ${\cal P}_{[n]}$,
denoted as ${\rm SA}_n$. 
We regard a generic random partition $\Pi_n\vdash [n]$ as a random allocation 
of balls labeled $1,\ldots,n$ to some set of nonempty boxes,
which the operation ${\rm SA}_n$ transforms into some other random allocation $\Pi_n'$.
Given $\Pi_n = \pi_n$,
\begin{itemize}
\item
if $K_n=1$, first delete a single ball picked uniformly at random from the balls
allocated according to $\pi_n$, 
to make an intermediate partition of some set of $n-1$ balls, 
then add to this intermediate partition a single box containing the deleted ball.

\item 
if $K_n=\{k_1,\ldots,k_r\}$, 
first pick out a sequence of $k_1-1$ of the $n$ balls from $\pi_n$ by 
uniform random sampling without replacement and 
put these $k_1-1$ balls together as set $\#1$; 
continue to pick out a sequence of $k_2-1$ of the remaining $n-k_1+1$ balls 
by uniform random sampling withour replacement and 
put these $k_2-1$ balls together as set $\#2$; 
keep doing these until we get $r$ set of balls marked as $\#1,\#2,\ldots,\#r$.

Then mark a ball chosen uniformly from the remaining 
$n-\sum_{i=1}^rk_i+r$ balls as ball $\#1$, 
continue to mark a ball chosen uniformly from remaining 
$n-\sum_{i=1}^rk_i+r-1$ unmarked balls as ball $\#2$, 
keep doing these until we marked $r$ balls. 
Now for each $i=1,2,\ldots,r$, 
add all balls in set $\#i$ into the box containing ball $\#i$.
\end{itemize}
In either case delete empty boxes in case any appear after the sampling
step.  
The resulting partition of $[n]$ is $\Pi_n'$. 
For each $q(n:\cdot)$, this defines a Markovian transition operator ${\rm SA}_n$ 
on partitions of $[n]$.

\begin{lemma}\label{generalstat} 
Let $\Pi_n$ be an exchangeable random partition of $[n]$ with finite
EPPF $p$ defined as a function of compositions of $m$ for $1 \le m \le n$.
Let $\Pi_n'$ be derived from $\Pi_n$ by the ${\rm SA}_n$ operation determined
by some sequence $q(n:\cdot)$ as {\rm \re{sageneralqb}}.
Then $\Pi_n'$ is an exchangeable random partition of $[n]$ whose EPPF
$p'$ is determined on compositions of $[n]$ by the formula
\begin{align}
\label{generalrec11}
&p'(n_1,n_2,\ldots,n_\ell)=
\frac{q(n:1)}{n}\sum_{j:n_j=1}\,p(\widehat{n_j})+\sum_{\{k_1,\ldots,k_r\}}q(n:k_1,\ldots,k_r;n-\sum_{j=1}^rk_j)\,\nonumber\\
& \sum_{\eta\in H_{\{k_1,\ldots,k_r\}}^{(n_1,\ldots,n_\ell)}}
\frac{\prod_{i=1}^\ell d(n_i; k_{\eta_{i(1)}}, \ldots, k_{\eta_{i(| \eta_i |)}};  n_i-\sum_{l=1}^{|\eta_i|} k_{\eta_{i(l)}} ) }{d(n; k_1,\ldots,k_r;n-\sum_{j=1}^rk_j)} p(n_1-\sum_{l=1}^{|\eta_1|} k_{\eta_{1(l)}} + |\eta_1|,\ldots, n_\ell-\sum_{l=1}^{|\eta_\ell|} k_{\eta_{\ell(l)}}+|\eta_\ell| )
\end{align}
where $(\widehat{n_j})$ is formed from $(n_1,n_2,\ldots,n_\ell)$ by deleting part $n_j$, 
$|\eta_i|$ is the number of elements in $\eta_i$ and 
the second sum on the right hand side is over all integer sets $\{k_1,\ldots,k_r\}$ 
with $r\ge 1, k_1,\ldots,k_r\ge 2, \sum_{i=1}^r k_i\le n$. 
\end{lemma}
\noindent
(Note that the right side of (\ref{generalrec11}) is identical to the right side of
\Mohle's recursion (\ref{generalrec1}).)

\proof
Let $K_n$ be a random element with its distribution according to the sequence $q(n:\cdot)$, 
i.e. $K_n$ equals to $1$ with probability $q(n:1)$, 
and equals to set $\{k_1,\ldots,k_r\}$ with probability $q(n:k_1,k_2,\ldots, k_r; s)$. 
For each partition $\pi_n'$ of $[n]$ we can compute
\begin{align}\label{generalrec111}
\prob ( \Pi_n' = \pi_n') &= q(n:1)\,\prob ( \Pi_n' = \pi_n' \giv K_n = 1) + \nonumber \\
&\sum_{\{k_1,\ldots,k_r\}} q(n:k_1,\ldots,k_r;n-\sum_{j=1}^rk_j) \,\prob ( \Pi_n' = \pi_n' \giv K_n = \{k_1,\ldots,k_r\}).
\end{align}
Assuming that $\pi_n'$ has boxes of sizes $n_1, \ldots, n_\ell$,
and that the ${\rm SA}_n$ operation acts on an exchangeable $\Pi_n$ with EPPF $p$, 
we deduce \re{generalrec11} from \re{generalrec111} and
\eq
\label{generalrec1one}
\prob ( \Pi_n' = \pi_n' \giv K_n = 1) = \frac{1}{n}\sum_{j:n_j=1} p(\ldots,\widehat{n_j},\ldots), 
\en
\begin{align}
\label{generalrec1111}
&\prob ( \Pi_n' = \pi_n' \giv K_n = \{k_1,\ldots,k_r\}) = \nonumber \\
&\sum_{\eta\in H_{\{k_1,\ldots,k_r\}}^{(n_1,\ldots,n_\ell)}}
\frac{\prod_{i=1}^\ell d(n_i; k_{\eta_{i(1)}}, \ldots, k_{\eta_{i(| \eta_i |)}};  n_i-\sum_{l=1}^{|\eta_i|} k_{\eta_{i(l)}} ) }{d(n; k_1,\ldots,k_r;n-\sum_{j=1}^rk_j)} p(n_1-\sum_{l=1}^{|\eta_1|} k_{\eta_{1(l)}} + |\eta_1|,\ldots, n_\ell-\sum_{l=1}^{|\eta_\ell|} k_{\eta_{\ell(l)}}+|\eta_\ell| ).
\end{align}
First consider \re{generalrec1111}.
By the definition, 
every $\eta$ in $H_{\{k_1,\ldots,k_r\}}^{(n_1,\ldots,n_\ell)}$ induces a class of ways to 
allocate sets of sampled balls into boxes of $\pi_n'$ with which 
the result $(\Pi_n' = \pi_n')$ is possible after the ${\rm SA}_n$ operation. 
With the particular class of allocation indicated by $\eta$, 
the actual sequence of labels of sampled and marked balls, in order of choice, 
can be any one of
$$
\prod_{i=1}^\ell \frac{n_i !}{(n_i-\sum_{l=1}^{|\eta_i|} k_{\eta_{i(l)}})! }
$$
sequences, out of
$$
\frac{n!}{(n-\sum_{j=1}^r k_j)!}
$$ 
total number of possibilities.

Given any particular order of sampling and marking balls, 
let $M_{\sum_{j=1}^r k_j - r}$ be the set of labels of the $\sum_{j=1}^r k_j - r$ balls that are moved.
Then the event $(\Pi_n' = \pi_n')$ occurs if and only if 
the restriction of $\Pi_n$ to $[n] -M_{\sum_{j=1}^r k_j - r}$ equals 
the restriction of $\pi_n'$ to $[n] - M_{\sum_{j=1}^r k_j - r}$, 
which is a particular partition of $n-\sum_{j=1}^r k_j + r$ labeled balls into boxes 
of $\bar{n}_1,  \ldots, \bar{n}_\ell$ balls, 
where $\bar{n}_i = n_i 1(|\eta_i|=0) + (n_i - \sum_{l=1}^{|\eta_i|}k_{\eta_i(l)} + |\eta_i|) 1(|\eta_i|\ne 0)$.
The conditional probability of $(\Pi_n' = \pi_n')$,
given $K_n$ equalling to some $\{k_1,\ldots,k_r\}$ and 
a particular order of sampling and marking balls with which 
the ${\rm SA}_n$ operation is performed, 
is therefore  
$$
p(n_1-\sum_{l=1}^{|\eta_1|} k_{\eta_{1(l)}} + |\eta_1|,\ldots, n_\ell-\sum_{l=1}^{|\eta_\ell|} k_{\eta_{\ell(l)}}+|\eta_\ell| )
$$
by the assumed exchangeability of $\Pi_n$, 
and the definition of the EPPF $p$ of $\Pi_n$ on compositions of $m \le n$ 
by restriction of $\Pi_n$ to subsets of size $m$.

Also by the definition of $H_{\{k_1,\ldots,k_r\}}^{(n_1,\ldots,n_\ell)}$, 
some choices $\eta$'s are counted once although they appear different here 
because of the labeling of elements in set $\{k_1,\ldots,k_r\}$, 
e.g. $(\{1,2\},\{3\})$, $(\{2,3\},\{1\})$ and $(\{1,3\},\{2\})$ are counted as one thing 
in $H_{\{3,3,3\}}^{(6,3)}$, 
but here since the balls are labeled according to the order they are sampled and marked, 
the repetition we eliminated in defining $H_{\{k_1,\ldots,k_r\}}^{(n_1,\ldots,n_\ell)}$ actually leads to 
different ways of sampling and adding, 
so we want to count these in by multiplying
$$
\prod_{j=2}^n {l_j \choose l_{\eta_1,j}, l_{\eta_2,j},\ldots,l_{\eta_\ell,j} }
$$ 
where $l_j=\#\{i: k_i=j\}$, $l_{\eta_m,j}:=\# \{i: k_{\eta_m(i)}=j \}$ for $j=2,\ldots,n$ 
and $m=1,\ldots,\ell$.

Now the evaluation  \re{generalrec1111} is apparent, and \re{generalrec1one} too is
apparent by a similar but easier argument.
\endpf

\par
The following proposition follows from above lemma and similar argument as in \cite{donggnedinpitman06}.

\begin{proposition}\label{generalsauniq}
For each sequence $q(n:\cdot)$ as (\ref{sageneralqb}), 
the corresponding ${\rm SA}_n$ transition operator on partitions of $[n]$ 
has a unique stationary distribution.
A random partition with this stationary distribution is
exchangeable, and its EPPF is the finite unique EPPF $p$ that
satisfies \Mohle's recursion {\rm (\ref{generalrec1})}, 
that is {\rm \re{generalrec11}} with $p'=p$.
\end{proposition}

\section{Infinite partitions}\label{generalINF}

Same as \cite[Section 6]{donggnedinpitman06}, 
in this section we pass from finite partitions to the projective limit, 
and arrive at the desired integral representation of infinite array $q^\infty$
satisfying recursion (\ref{generalqrec1}), (\ref{generalqrec1b}). 
We get the main result which is the infinite counterpart of Theorem \ref{finitegeneralrec}, 
and can be seen as generalized version of \cite[Theorem 6.2]{donggnedinpitman06}.

\par 
An infinite sequence of freeze-and-merge operations $\FM:=(\FM_n, n=1,2,\ldots)$ 
which satisfies the condition in Definition \ref{defqconsis} for 
all positive integers $1\le m<n<\infty$ is called {\em consistent}. 
For each $n=1,2,\ldots$ the Markov chain starting from $\Sigma_n^*$ and 
driven by $\FM_n$ terminates with an induced final partition $\Pi_n$.
These comprise an infinite partition $\Pi_\infty=(\Pi_n)_{n=1}^\infty$ which 
we call the {\em final partition} associated with consistent infinite $\FM$.

\begin{lemma}\label{infinitexi} 
For every infinite array 
\eq\label{generalqsequenceinf}
q^{\infty}:=(q(1:\cdot),q(2:\cdot),\ldots,q(n:\cdot),\ldots)
\en
where 
\eq
q(1:\cdot)=\{q(1:1)\}=\{1\}
\en
and for $b>1$, $q(b:\cdot)$ are as in {\rm \re{generalqb}}, 
with entries satisfying the recursion {\rm(\ref{generalqrec1}), (\ref{generalqrec1b})}, 
there exist a non-negative finite measure on the infinite simplex with the form 
$\Xi=\Xi_0+a\delta_0$, 
where $\Xi_0$ has no atom at zero and $\delta_0$ is a unit mass at zero, 
and a non-negative real number $\rho$ such that the entries of $q^\infty$ can 
be represented by $(\Xi,\rho)$ as {\rm \re{xiqsequencek}, \re{xiqsequence1}}. 
The data $(\Xi,\rho)$ are unique up to a positive factor.
\end{lemma}

\proof
Suppose $q$ solves the recursion (\ref{generalqrec1}), (\ref{generalqrec1b}), 
and suppose $q(2:2;0)<1$.
Let $\Phi(n), n=1,2,\ldots$ satisfy
\eq\label{infinitexipf1}
{\Phi(n)\over\Phi(n+1)}=1-{1\over n+1}q(n+1:1)-{2\over n+1}q(n+1:2;n-1)
\en
for $n\geq 1$; since the right hand side is strictly positive this recursion has a unique  
solution with some given initial value $\Phi(1)=\rho$, where $\rho>0$.
For each $q(n; k_1,k_2,\ldots,k_r; s)$ set
$$
\lambda_{n; k_1,k_2,\ldots,k_r; s}:=\frac{q(n:k_1,k_2,\ldots,k_r; s)}{d(n;k_1,k_2,\ldots, k_r; s)}\Phi(n)
$$ 
then from (\ref{infinitexipf1}) and (\ref{generalqrec1}), we can derive \re{xilambdaconsis}: 
$$
\lambda_{n; k_1,k_2,\ldots,k_r; s}=\sum_{i=1}^r \lambda_{n; k_1,\ldots,k_{i-1},k_i+1,k_{i+1},\ldots,k_r; s}
+s\lambda_{n+1; k_1,k_2,\ldots,k_r,2; s-1}+\lambda_{n+1; k_1,k_2,\ldots,k_r; s+1}
$$
then by \cite[Lemma 18]{jason00} and \cite[Theorem 2]{jason00} 
we conclude \re{xilambdaint}, hence \re{xiqsequencek} holds for some non-negative finite measure 
on the infinite simplex with the form $\Xi=\Xi_0+a\delta_0$, 
where $\Xi_0$ has no atom at zero and $\delta_0$ is a unit mass at zero, 

From (\ref{generalqrec1b}) we find
$$
\rho={\Phi(1)q(1:1)\over 1 }=\cdots={\Phi(n)q(n:1)\over n}=\cdots,
$$
hence by setting $\Phi(n:1):=\rho n$ we deduce \re{xiqsequence1}.
For the special case $q(2:2;0)=1$, it is easy to observe that $\rho =0$, 
and we get $\Xi=\delta_0$ by similar analysis.
\endpf
\vskip 0.5cm

\par Recording this lemma together with previous results, we have the following result, 
which is the counterpart of \cite[Theorem 6.2]{donggnedinpitman06}:

\begin{theorem} \label{generalmain}
Let $\Pi_\infty=(\Pi_n)_{n=1}^\infty$ be a nontrivial exchangeable random partition of $\Nat$, 
different from the trivial one-block partition.
The following are equivalent:
\begin{itemize}
\item[{\rm(i)}] The EPPF $p$ satisfies M\"ohle's recursion {\rm (\ref{generalrec1})} with 
some infinite array $q^\infty$ with form {\rm\re{generalqsequenceinf}}. 
\item[{\rm(ii)}] This array is representable as {\rm \re{xiqsequencek}, \re{xiqsequence1}} 
by some nontrivial $(\Xi,\rho)$ which is unique up to a positive factor, 
as claimed in Lemma {\rm \ref{infinitexi}}.
\item [{\rm(ii)}] This $\Pi_\infty$ is induced by the final partition of some standard 
$\Xi$-coalescent freezing at rate $\rho$.
\item[{\rm (iii)}] This $\Pi_\infty$ is the final partition of some consistent infinite $\FM$ operation.
 \end{itemize}
\end{theorem}

\vskip 0.5cm

\par Finally, we complete this paper with the following uniqueness assertion, 
similar with \cite[Lemma 6.3]{donggnedinpitman06}:

\begin{lemma}\label{generalinj}
The correspondence $q^\infty \mapsto p$ between infinite arrays with 
$q(2:1)>0$ satisfying consistency {\rm  (\ref{generalqrec1}), (\ref{generalqrec1b})} 
and the EPPF's is bijective.  
\end{lemma}
\proof
We only need to show that $p$ uniquely determines $q$. 
For general infinite partitions, $q(2:1)=p(1,1)>0$ implies that $p(1,1,\ldots,1)>0$. 
By Lemma \ref{fmfinal}, $p$ must solve M\"ohle's recursion \re{generalrec1}, so 
$$p(1,\ldots,1)=q(n:1)q(n-1:1)\cdots q(2:1)$$
shows that the $q(n:1)$'s are uniquely determined by $p$. 
By exploiting the formula 
\begin{align}
p(m,1,\ldots,1)=&{q(n:m;n-m)\over{n\choose m}}p(1,\ldots,1)+q(n:1)\frac{n-m}{n}p(m,\widehat{1},1,\ldots,1)\nonumber\\
&+\sum_{k=2}^{m-1}q(n:k;n-k)\frac{{m\choose k}}{{n\choose k}}p(m-k+1, 1,\ldots,1)
\end{align}
with induction in $m=2,3,\ldots,n-1$, 
it is clear that entries with form $q(n:m;n-m)$, $2\le m\le n-1$ are also uniquely determined by $p$. 
Similarly, if we look at the equation \re{generalrec1} with $p(m,l,1,\ldots,1)$, $2\le m,\ l\le n$, $m+l\le n$ 
on the left hand side, 
by induction in $m$, $l$, 
we can deduce that entries with form $q(n: m,l;n-m-l)$ are uniquely determined by $p$ as well. 
Same mechanism can be carried on to conclude that all entries of $q^\infty$ are uniquely determined by $p$.

The uniqueness fails when $q(2:1)=0$, which corresponds to the case with no freezing, 
in that case the singleton partition will be the final partition regardless of the coalescing theme. 
\endpf

\vskip 0.5cm

\section*{Acknowledgments}  
Most of this work was carried out during the progress of \cite{donggnedinpitman06}. 
The author is very grateful to Alexander Gnedin and Jim Pitman for insights and 
many valuable discussions.

\end{document}